\documentclass[12pt]{amsart}
\usepackage{amssymb}
\usepackage[latin1]{inputenc}
\usepackage{amsthm}
\usepackage{amsfonts}
\usepackage{mathrsfs}
\usepackage{graphicx}
\usepackage{amsmath}
\usepackage[all]{xy}

\let\mathcal\mathscr

\def\bC{{\mathbb C}}

\def\bP{{\mathbb P}}
\def\bN{{\mathbb N}}
\newtheorem{thm}{Theorem}[section]

\def\tilde{\widetilde}

\def\phi{\varphi}

\numberwithin{equation}{section}

\newtheorem{theorem}[thm]{Theorem}

\newtheorem{conjecture}[thm]{Conjecture}
\newtheorem{corollary}[thm]{Corollary}

\newtheorem{definition}[thm]{Definition}
\newtheorem{example}[thm]{Example}

\newtheorem{lemma}[thm]{Lemma}

\newtheorem{proposition}[thm]{Proposition}

\newtheorem{remark}[thm]{Remark}

\begin{document}


\title {Hyperbolicity of geometric orbifolds}
\author{Erwan Rousseau}
\date{}

\begin{abstract}
We study complex hyperbolicity in the setting of geometric orbifolds introduced by F. Campana. Generalizing classical methods to this context, we obtain degeneracy statements for entire curves with ramification in situations where no Second Main Theorem is known from value distribution theory.
\end{abstract}
\maketitle

\tableofcontents

\section{Introduction}
F. Campana has introduced in \cite{C04} orbifold structures, namely pairs $(X/\Delta)$ with $X$ a complex manifold and a divisor $\Delta=\sum_{i}(1-\frac{1}{m_{i}})Z_{i}$ where the $Z_{i}$ are distinct irreducible divisors and $m_{i} \in \bN\cup \{\infty\}$, as a new frame for the classification of compact K\"ahler manifolds. 
These structures appeared naturally for fibrations $f:X\rightarrow Y$. Indeed the multiple fibres of $f$ lead to the definition of the orbifold base of $f$, $(Y/\Delta(f))$ where
$$\Delta(f):=\sum_{D \subset Y} \left(1-\frac{1}{m(f,D)}\right)D$$
$m(f,D)$ being the multiplicity of the fiber of $f$ above the generic point of $D$.
 A new class of varieties was then introduced, the \textit{special varieties}, as the varieties which do not admit fibrations of general type i.e with an orbifold base of general type. Campana \cite{C04} proves the existence for every complex algebraic manifold $X$ of a fibration $c_{X}:X \rightarrow C(X)$, the core of $X$, such that its general fibers are special and if $X$ is not special, $c_{X}$ is of general type.

These geometric orbifolds should be considered as true geometric objects as one can define for them differential forms, fundamental groups, Kobayashi pseudo-distance... Here we study the hyperbolic aspects of these objects. An important conjecture of Campana \cite{C04} is that $X$ is special if and only if the Kobayashi pseudo-distance $d_{X}$ vanishes identically on $X \times X$. This is known only for curves, projective surfaces not of general type and rationally connected manifolds. 

This conjecture then implies that $d_{X}$ should be the pull-back by $c_{X}$ of the Kobayashi pseudo-distance $\delta_{X}$ of the orbifold base of the core.

The study of the hyperbolic aspects of one-dimensional orbifolds has been done in \cite{CW}. Here we study hyperbolicity of higher dimensional orbifolds following the philosophy of Campana that one should study these objects generalizing the tools we use for manifolds without orbifold structures or logarithmic manifolds. The paper is organized as follows.

In section 2, we recall the basic facts on geometric orbifolds following \cite{C04} and \cite{C07}. 

In section 3, we recall the definitions of classical and non-classical Kobayashi hyperbolicity for orbifolds. Then we illustrate these notions in the case of orbifold curves. We compute explicitly the orbifold Kobyashi pseudo-distance for $$(X/\Delta)=(D/(1-\frac{1}{n})\{0\}), 0<n\in  \bN\cup \{\infty\},$$ where $D$ is the unit disk. This answers a question of Campana and Winkelmann (see \cite{CW}) and enables us to recover as a corollary the equivalence of classical and non-classical hyperbolicity for orbifold curves. Finally, we show that this is not the case in higher dimension giving an example of an orbifold surface which is classically hyperbolic but not hyperbolic.

In section 4, we define and study algebraic hyperbolicity in the orbifold setting. We prove that $(\mathbb{P}^n/\Delta)$ is algebraically hyperbolic where $\Delta = \sum_{1 \leq i \leq q} (1-\frac{1}{m_{i}})H_{i}$ for $H_{1},\dots,H_{q}$ very generic hypersurfaces and $\deg \Delta >2n$.

In section 5, we discuss an orbifold Kobayashi's conjecture motivated by the results of the preceding section.

In section 6, we define and use orbifold jet differentials. The main applications are algebraic degeneracy statements for entire curves with ramification in situations where no Second Main Theorem is known from value distribution theory. Namely, we prove
\begin{theorem}\label{theo1}
Let $(X/\Delta)$ be a smooth projective orbifold surface of general type where $\Delta$ has the following decomposition into irreducible components, $\Delta=\sum_{i=1}^{n}(1-\frac{1}{m_{i}})C_{i}$. Suppose that $g_{i}:=g(C_{i}) \geq 2$, $h^0(C_{i},\mathcal{O}_{C_{i}}(C_{i})) \neq 0$ for all $i$ and
that the logarithmic Chern classes of $(X, \lceil \Delta \rceil)$ verify
\begin{equation}
\overline{c_{1}}^2-\overline{c_{2}}-\sum_{i=1}^{n}\frac{1}{m_{i}}(2g_{i}-2+\sum_{j \neq i}C_{i}C_{j})>0,
\end{equation}
then there exists a proper subvariety $Y\subsetneq X$ such that every entire curve $f: \bC \rightarrow X$ which is an orbifold morphism, i.e ramified over $C_{i}$ with multiplicity at least $m_{i}$, verifies $f(\bC) \subset Y$.
\end{theorem}

This result can be seen as an orbifold version of results of McQuillan \cite{McQ0} (see also \cite{McQ2} and \cite{EG} for the logarithmic case) on the Green-Griffiths-Lang conjecture which can be generalized to the orbifold setting
\begin{conjecture}
Let $(X/\Delta)$ be a smooth projective orbifold of general type. Then there exists a proper subvariety $Y\subsetneq X$ such that every orbifold morphism $f: \bC \rightarrow (X/\Delta)$ verifies $f(\bC) \subset Y$.
\end{conjecture}

The methods used also enable us to generalize a result of Campana and Paun on weakly-special manifolds \cite{CP}.

In section 7, we study measure hyperbolicity of orbifolds.

\bigskip

\textbf{Acknowledgements}. We would like to thank Fr\'ed\'eric Campana and Michael McQuillan for their suggestions and their interests in this work.

\section{Geometric orbifolds, special varieties and hyperbolicity}
\subsection{Geometric orbifolds}
Let $X$ be a complex manifold. Following F. Campana (see \cite{C04}, \cite{C07}) we recall the basic facts on orbifolds and why these structures should be useful to describe the hyperbolic aspects of projective manifolds.
\begin{definition}
A (geometric) orbifold $(X/\Delta)$ is a pair consisting of $X$ with a $\mathbb{Q}$-divisor $\Delta$ on $X$ for which the decomposition in reduced irreducible divisors is of the form
\begin{displaymath}
\Delta=\sum_{i} \left(1-\frac{1}{m_{i}}\right)Z_{i},
\end{displaymath}
where $0<m_{i} \in \bN\cup \{\infty\}$.
\end{definition}

These structures appear naturally in the case of fibrations.

\begin{definition}
Let $f: X \rightarrow Y$ be a fibration between complex compact manifolds. For every irreducible divisor $D \subset Y$:
 $$f^{*}(D)=\sum_{j}m_{j}D_{j}+R$$
where $R$ denotes the $f$-exceptional part. Then the multiplicity of $D$ is defined by
$$m(f,D)=\inf_{j}\{m_{j}\}.$$
\end{definition}

\begin{definition}
The orbifold base $(Y/\Delta(f))$ is then defined by
$$\Delta(f)=\sum_{D \subset Y}\left(1-\frac{1}{m(f,D)}\right)D.$$
\end{definition}

\begin{remark}
This can also be defined for meromorphic fibrations, resolving indeterminacies. As we will work up to bimeromorphic equivalence, we will not specify when a fibration is meromorphic.
\end{remark}

\begin{definition}
\begin{enumerate}

\item The canonical line bundle of $(X/\Delta)$ is defined by $K_{(X/\Delta)}:=K_{X}+\Delta.$
\item The Kodaira dimension of $(X/\Delta)$ is $\kappa(X/\Delta):=\kappa(K_{(X/\Delta)})$
\item Let $f: X \rightarrow Y$ be a fibration between complex compact manifolds. Then its Kodaira dimension is $$\kappa(f)=\inf_{f'\sim f}\{\kappa(Y'/\Delta(f'))\},$$
where $f'\sim f$ iff there exists a commutative diagramm
$$
\xymatrix{
    X' \ar[r]^w \ar[d]_{f'}& X \ar[d]^{f} \\
    Y' \ar[r]_v & Y
  }
$$
with $w$ and $v$ bimeromorphic.
\end{enumerate}
\end{definition}

Now, one can introduce the \textit{special} geometry. 

\subsection{Special varieties and hyperbolicity}
Let $X$, $Y$ be projective manifolds.

\begin{definition}
\begin{enumerate}
\item A fibration $f: X \rightarrow Y$ is of general type if $\kappa(f)=\dim Y>0.$

\item $X$ is \textit{special} if there is no fibration $f: X \rightarrow Y$ of general type.

\item $f: X \rightarrow Y$ is special if its general fiber is special.
\end{enumerate}
\end{definition}

\begin{example}
Rationnally connected manifolds and manifolds with Kodaira dimension $0$ are two important examples of special varieties (see \cite{C04}).
\end{example}

The following theorem makes clear why orbifold structures are useful for hyperbolicity.
\begin{theorem}[Campana \cite{C04}]

Let $X$ be a projective manifold. Then there exists a unique (up to equivalence) fibration 
$$c_{X}:X \rightarrow \mathcal{C}(X)$$ called the \textit{core} of $X$
such that 
\begin{enumerate}

\item $c_{X}$ is special.

\item $c_{X}$ is of general type or constant (iff $X$ is special).
\end{enumerate}
\end{theorem}

In other words, one can decompose a projective manifold into its "hyperbolic" part, the orbifold base of the core, and its "non-hyperbolic" part. Indeed, special varieties should be characterized by the following conjecture:
\begin{conjecture}[Campana \cite{C04}]
$X$ is special iff $d_{X} \equiv 0,$ where $d_{X}$ denotes the Kobayashi pseudo-distance.
\end{conjecture}

\begin{remark}
It is only known to be true for curves, projective surfaces not of general type, rationally connected manifolds.
\end{remark}

As a consequence, one should have the following description of the Kobayashi pseudo-distance:
\begin{conjecture}[Campana \cite{C04}]\label{Kob0}
$d_{X}=c_{X}^*\delta_{X}$ where $\delta_{X}$ is a pseudo-distance on $\mathcal{C}(X)$.
\end{conjecture}
We will give some precisions to what that pseudo-distance $\delta_{X}$ should be in the next section.

\section{Kobayashi hyperbolicity of orbifolds}
\subsection{Orbifold Kobayashi pseudo-distance}
First, let us recall following \cite{CW} the definition of classical and non-classical orbifold morphisms from the unit disk to an orbifold.
\begin{definition} 
Let $(X/\Delta)$ be an orbifold with $\Delta=\sum_{i} (1-\frac{1}{m_{i}})Z_{i}$, $D=\{z \in \bC / |z|<1\}$ the unit disk and $h$ a holomorphic map from $D$ to $X$.
\begin{enumerate}
\item $h$ is a (non-classical) orbifold morphism from $D$ to $(X/\Delta)$ if $h(D) \nsubseteq supp(\Delta)$ and $mult_{x}(h^*Z_{i})\geqslant m_{i}$ for all $i$ and $x\in D$ with $h(x)\in supp(Z_{i})$. If $m_{i}=\infty$ we require $h(D)\cap Z_{i} = \emptyset$.

\item $h$ is a classical orbifold morphism from $D$ to $(X/\Delta)$ if $h(D) \nsubseteq supp(\Delta)$ and $mult_{x}(h^*Z_{i})$ is a multiple of $m_{i}$ for all $i$ and $x\in D$ with $h(x)\in supp(Z_{i})$. If $m_{i}=\infty$ we require $h(D)\cap Z_{i} = \emptyset$.

\end{enumerate}
\end{definition}

Then we can define classical and non-classical orbifold morphisms between orbifolds.

\begin{definition}
Let $(X/\Delta)$ and  $(X'/\Delta')$ be orbifolds and $f: X \rightarrow X'$ a holomorphic map. $f$ is an orbifold morphism (resp. classical orbifold morphism) from $(X/\Delta)$ to  $(X'/\Delta')$ if 
\begin{enumerate}
\item $f(X) \nsubseteq supp(\Delta')$.
\item $f\circ g : D \rightarrow (X'/\Delta')$ is an orbifold morphism (resp. classical orbifold morphism) for all orbifold morphism (resp. classical orbifold morphism) $g: D \rightarrow (X/\Delta)$ such that $g(D)\nsubseteq f^{-1}(supp(\Delta'))$.
\end{enumerate}
\end{definition}

Let $(X/\Delta)$ be an orbifold with $\Delta=\sum_{i} a_{i}Z_{i}$ and $\Delta_{1}$ the union of all $Z_{i}$ with $a_{i}=1$. 
\begin{definition}
\begin{enumerate}
\item The orbifold Kobayashi pseudo-distance $d_{(X/\Delta)}$ on $(X/\Delta)$ is the largest pseudo-distance on $X\setminus \Delta_{1}$ such that
\begin{displaymath}
g^*d_{(X/\Delta)} \leqslant d_{P}
\end{displaymath}
for every orbifold morphism $g: D \rightarrow (X/\Delta)$, where $d_{P}$ denotes the Poincar\'e distance on $D$.
\item The classical orbifold Kobayashi pseudo-distance $d^*_{(X/\Delta)}$ on $(X/\Delta)$ is the largest pseudo-distance on $X\setminus \Delta_{1}$ such that
\begin{displaymath}
g^*d_{(X/\Delta)} \leqslant d_{P}
\end{displaymath}
for every classical orbifold morphism $g: D \rightarrow (X/\Delta)$, where $d_{P}$ denotes the Poincar\'e distance on $D$.
\end{enumerate}
\end{definition}

As an immediate consequence of the definition we have
\begin{proposition}\label{dec}
Let $f:(X/\Delta) \rightarrow  (X'/\Delta')$ be an orbifold morphism (resp. classical orbifold morphism). Then $f^*d_{(X'/\Delta')}\leq d_{(X/\Delta)}$ (resp. $f^*d^*_{(X'/\Delta')}\leq d^*_{(X/\Delta)}$).
\end{proposition}

\begin{definition}
An orbifold  $(X/\Delta)$ is hyperbolic (resp. classically hyperbolic) if $d_{(X/\Delta)}$ (resp.  $d^*_{(X/\Delta)}$) is a distance on $X\setminus \Delta_{1}$. 
\end{definition}

A corollary of proposition \ref{dec} is
\begin{corollary}\label{brody1}
Let $(X/\Delta)$ be a hyperbolic (resp. classically hyperbolic) orbifold. Then every orbifold morphism (resp. classical orbifold morphism) $f: \bC \rightarrow (X/\Delta)$ is constant.
\end{corollary}

In the compact and logarithmic setting, Brody-type theorem turn out to be very useful to characterize hyperbolicity as converse of corollary \ref{brody1}. This is done in \cite{CW} with the following theorem
\begin{theorem}\label{brody2}
Let $(X/\Delta)$ be a non hyperbolic (resp. non classically hyperbolic) compact orbifold. Then there exists either a non constant orbifold morphism (resp. classical orbifold morphism) $f: \bC \rightarrow (X/\Delta)$ or a non-constant holomorphic map $f: \bC \rightarrow supp(\Delta)$.
\end{theorem}

Then we can refine conjecture \ref{Kob0}
\begin{conjecture}[Campana \cite{C04}]\label{Kob}
$d_{X}=c_{X}^*\delta_{X}$ where $\delta_{X}=d_{(\mathcal{C}(X)/\Delta(c_{X}))}$.
\end{conjecture}
\subsection{Hyperbolicity of orbifold curves}
Let us illustrate orbifold hyperbolicity in dimension 1. This has been studied in \cite{CW}. We give here a different approach. The following result gives a concrete example where one can compute the (non-classical) orbifold Kobayashi pseudo-distance, answering a question of \cite{CW}.

\begin{theorem}\label{curv}
Let $(X/\Delta)=(D/(1-\frac{1}{n})\{0\})$, $0<n\in  \bN\cup \{\infty\}$. Then
$$d_{(X/\Delta)}=d^*_{(X/\Delta)}$$
obtained as the pseudo-distance $d$ induced by $$\omega=\frac{4dzd\overline z}{n^2|z|^{2-\frac{2}{n}}(1-|z|^{\frac{2}{n}})^2}.$$
\end{theorem}
\begin{remark}
For $n=\infty$ we obtain $\omega=\frac{4dzd\overline z}{|z|^{2}(\log|z|^{2})^2}$, which is known to induce the Kobayashi distance on the punctured disc.
\end{remark}
To prove this we shall need the Ahlfors-Schwarz lemma (see e.g \cite{De95})
\begin{lemma}[Ahlfors-Schwarz lemma]
Let $\gamma(t)=\gamma_{0}(t)dt d\overline t$ be a singular hermitian metric on $D$ where $\log \gamma_{0}$ is a subharmonic function such that $i\partial\overline{\partial}\log \gamma_{0}(t)\geq Ai\gamma_{0}(t)dt\wedge d\overline t$ in the sense of currents, for some positive constant $A$. Then
$$\gamma(t)\leq \frac{2}{A}\frac{dzd\overline z}{(1-|z|^2)^2}=\frac{1}{2A}h_{P}$$
where $h_{P}$ denotes the Poincar\'e metric.
\end{lemma}
Now we prove the claim of the previous theorem \ref{curv}
\begin{proof}
As $z\rightarrow z^n$ gives an unfolding $D \rightarrow (X/\Delta)$ (see \cite{CW}), the classical Kobayashi pseudo-distance is the push-forward of the Poincar\'e metric on $D$ which is $\frac{4dzd\overline z}{n^2|z|^{2-\frac{2}{n}}(1-|z|^{\frac{2}{n}})^2}.$ Therefore we have $d_{(X/\Delta)}\leq d.$

Now, let $f:D \rightarrow (X/\Delta)$ be an orbifold morphism. We shall apply the Ahlfors-Schwarz lemma to $f^*\omega$. We have $$f^*\omega=\gamma_{0}(t)dt d\overline t:=\frac{4|f'(t)|^2}{n^2|f(t)|^{2-\frac{2}{n}}(1-|f(t)|^{\frac{2}{n}})^2}dt d\overline t.$$
We remark that $\left(\frac{f'(t)}{f(t)^{1-\frac{1}{n}}}\right)^n$ is holomorphic. Indeed take $t_{0}$ such that $f(t_{0})=0$ and take local coordinate $t$ centered at $t_{0}$. Then $f(t)=t^mg(t)$ with $m\geq n$ and $g(0)\neq 0$. Therefore $f'(t)=t^{m-1}h(t)$ and $$\left(\frac{f'(t)}{f(t)^{1-\frac{1}{n}}}\right)^n=\frac{t^{n(m-1)}h^n(t)}{t^{m(n-1)}g^{n-1}(t)}=t^{m-n}\frac{h^n(t)}{g^{n-1}(t)}.$$
So $f^*\omega$ is a singular hermitian metric with $\log \gamma_{0}$ subharmonic.
Outside $f^{-1}(0)$ we have $$i\partial\overline{\partial}\log \gamma_{0}(t)=f^*(-Ricci(\omega))=f^*(\frac{1}{2}\omega).$$
Therefore by Ahlfors-Schwarz lemma
$$f^*\omega \leq  h_{P}$$
and finally  $d_{(X/\Delta)}\geq d.$

\end{proof}

Let us recall that the uniformization of orbifold curves is well known. 

\begin{theorem}[see \cite{Far}, p.234]\label{unif}
Let $(M/\Delta)$ be a compact orbifold curve with $$\Delta:=\sum_{i} (1-\frac{1}{\nu_{i}})x_{i}.$$

If $M=\mathbb{P}^1 $ we exclude two cases:

\begin{enumerate}
\item[(i)] $\{x_{1},x_{2},\dots\}$ consists of one point and $\nu
_{1}\neq \infty .$ \item[(ii)] $\{x_{1},x_{2},...\}$ consists of
two points and $\nu _{1}\neq \nu _{2}.$
\end{enumerate}

Let $M^{\prime }=M\setminus \Delta_{1}.$ Then there exists a simply connected Riemann surface
$\widetilde{M},$ a Kleinian group $G$ of self mappings of
$\widetilde{M}$ such that
\begin{enumerate}
\item[(a)] $\widetilde{M}/G\cong M^{\prime }$

\item[(b)] the natural projection $\pi :\widetilde{M}\rightarrow
M^{\prime }$ is
unramified except over the points $x_{k}$ with $\nu _{k}<\infty $ where the branch numbers verify $%
b_{\pi }(\widetilde{x})=\nu _{k}-1$ for all $\widetilde{x}\in \pi
^{-1}(\{x_{k}\}).$
\end{enumerate}
Moreover,
\begin{enumerate}
\item $\widetilde{M}\cong \bP^1$ iff $\deg(K_{(M/\Delta)})<0$
\item $\widetilde{M}\cong \bC$ iff $\deg(K_{(M/\Delta)})=0$
\item $\widetilde{M}\cong D$ iff $\deg(K_{(M/\Delta)})>0$.
\end{enumerate}
\end{theorem}

As a corollary we obtain
\begin{corollary}
Let $(X/\Delta)$ be a classically hyperbolic compact orbifold curve. Then the classical and non-classical Kobayashi pseudo-distances coincide and therefore $(X/\Delta)$ is hyperbolic.
\end{corollary}
\begin{proof}
By theorem \ref{unif}, $(X/\Delta)$ is uniformized by a Galois covering $D\rightarrow (X/\Delta)$ which ramifies exactly on $\Delta$. So the classical Kobayashi pseudo-distance is induced by $\omega$ the push-forward of the Poncar\'e metric. Let $f:D \rightarrow (X/\Delta)$ be an orbifold morphism. As the situation is locally the same as in theorem \ref{curv}, we obtain that $f^*\omega \leq  h_{P}$ and therefore $d_{(X/\Delta)}\geq d^*_{(X/\Delta)}$.
\end{proof}

Let us underline another property of hyperbolic orbifold curves which will be useful in the next section. In the case where $(M/\Delta)$ is uniformized by the unit disk $D$, we can project the Poincar\'{e} metric $\frac{4\left| dz\right| ^{2}}{%
(1-\left| z\right| ^{2})^{2}}$ on $M^{\prime }=M\setminus \Delta_{1}$ which gives a singular metric $h_{(M/\Delta)}$ inducing the Kobayashi pseudodistance. Despite the singularities we have

\begin{theorem}[\cite{Far}, p.233]\label{area}The area of $M^{\prime }=M\setminus \Delta_{1}$ with respect to the metric  $h_{(M/\Delta)}$ is finite and
\begin{equation*}
Area(M^{\prime})=2\pi \deg(K_{(M/\Delta)}).
\end{equation*}
\end{theorem}

 We have obtained a different proof of the result of \cite{CW} where the authors prove the equivalence of classical and non-classical hyperbolicity for orbifolds of dimension 1 using the orbifold Brody's lemma and Nevanlinna theory:
\begin{theorem}[\cite{CW}]\label{hypcur}
Let $(X/\Delta)$ be a compact orbifold of dimension 1. Then the following properties are equivalent
\begin{enumerate}
\item $(X/\Delta)$ is hyperbolic.
\item $(X/\Delta)$ is classically hyperbolic.
\item $\deg(K_{(X/\Delta)}):=\deg(K_{X}+\Delta)>0.$
\end{enumerate}
\end{theorem}

As a particular case, we recover the following classical result
\begin{theorem}[Nevanlinna \cite{Nev}]
Let $f:\bC \rightarrow (\bP^1/\sum_{i}(1-\frac{1}{m_{i}})a_{i})$ be a non-constant orbifold morphism. Then
$$\sum_{i} \left(1-\frac{1}{m_{i}}\right) \leq 2.$$
\end{theorem}

This theorem can be obtained thanks to the Second Main Theorem on $\bP^1$. In section 6 of this paper, we will obtain similar results in higher dimension in cases where no Second Main Theorem is known.

\subsection{Hyperbolicity of higher dimensional orbifolds}
Now we investigate the higher dimensional case. In particular,  is it still true that classical hyperbolicity coincides with (non-classical) hyperbolicity ?
The following example found in discussion with F. Campana answers it in the negative:

\begin{theorem}
There exists an orbifold surface $(S/\Delta)$ which is classically hyperbolic but not hyperbolic.
\end{theorem}

\begin{proof}
Let $X$ be a hyperbolic projective complex surface, $S$ its blow-up at one point $p$ and $E \subset S$ the exceptional divisor. We take $C_{1}, C_{2}, C_{3}$ three curves tangent to $E$ at three distinct points $p_{1}, p_{2}, p_{3}$. This can be done taking the strict transforms of curves in $X$ with a cusp at $p$ of local equation $(y+tx)^2-x^3=0$ for $t=0,1,-1.$ Define 
\begin{displaymath}
\Delta=(1-\frac{1}{3})C_{1}+(1-\frac{1}{3})C_{2}+(1-\frac{1}{5})C_{3}.
\end{displaymath}
Let's prove that  $(S/\Delta)$ is classically hyperbolic but not hyperbolic. 

$(E/\Delta')$, with 
\begin{displaymath}
\Delta'=(1-\frac{1}{2})p_{1}+(1-\frac{1}{2})p_{2}+(1-\frac{1}{3})p_{3},
\end{displaymath}
is not hyperbolic by theorem \ref{hypcur} since 
\begin{displaymath}
-2+(1-\frac{1}{2})+(1-\frac{1}{2})+(1-\frac{1}{3})=-\frac{1}{3}<0. 
\end{displaymath}
A non-constant orbifold morphism $f: \bC \rightarrow (E/\Delta')$ gives a non-constant orbifold morphism $f: \bC \rightarrow (S/\Delta)$ since the multiplicities verify respectively $2\times2\geq 3, 2\times2\geq 3, 2\times3\geq 5$. Therefore $(S/\Delta)$ is not hyperbolic.

Suppose $(S/\Delta)$ is not classically hyperbolic. Then by theorem \ref{brody2} there is either a non constant classical orbifold morphism $f: \bC \rightarrow (S/\Delta)$ or a non-constant holomorphic map $f: \bC \rightarrow supp(\Delta)$. But since $X$ is hyperbolic we must have $f(\bC) \subset E$. Therefore we obtain a non-constant classical orbifold morphism $f: \bC \rightarrow (E/\Delta')$ with 
\begin{displaymath}
\Delta'=(1-\frac{1}{m_{1}})p_{1}+(1-\frac{1}{m_{2}})p_{2}+(1-\frac{1}{m_{3}})p_{3}
\end{displaymath}
a divisor on $E$ such that the multiplicities verify, by the definition of a classical morphism, that $2m_{1}$ and $2m_{2}$ are multiples of 3, $2m_{3}$ is a multiple of 5. So, $m_{1}$ and $m_{2}$ are multiples of 3, $m_{3}$ is a multiple of 5 and we obtain
\begin{displaymath}
-2+(1-\frac{1}{m_{1}})+(1-\frac{1}{m_{2}})+(1-\frac{1}{m_{3}}) \geq -2+(1-\frac{1}{3})+(1-\frac{1}{3})+(1-\frac{1}{5})=\frac{2}{15}>0.
\end{displaymath}
This implies that $(E/\Delta')$ is hyperbolic and $f: \bC \rightarrow (E/\Delta')$ should be constant, which is a contradiction. So, $(S/\Delta)$ is classically hyperbolic.
\end{proof}

Now, we would like to define hyperbolic imbedding for orbifolds.

\begin{definition}
We say that $(X/\Delta)$ is hyperbolically (resp. classically hyperbolically) imbedded in $X$ if for any two sequences of points $(p_{n}), (q_{n}) \subset X\setminus \Delta_{1}$ converging to two points $p,q \in X$
\begin{displaymath}
d_{(X/\Delta)}(p_{n},q_{n}) \underset{n \rightarrow +\infty} \longrightarrow 0 \Rightarrow p=q
\end{displaymath}
(resp. $d^*_{(X/\Delta)}(p_{n},q_{n}) \underset{n \rightarrow +\infty} \longrightarrow 0 \Rightarrow p=q).$
\end{definition}

This can be characterized by
\begin{proposition}
Let $\omega$ be a hermitian metric on $X$ compact. Then $(X/\Delta)$ is hyperbolically (resp. classically hyperbolically) imbedded in $X$ iff there is a positive constant $c$ such that
\begin{displaymath}
f^*\omega\leq c\ h_{P}
\end{displaymath}
for all orbifold (resp. classical orbifold) morphism $f:D\rightarrow (X/\Delta)$, where $h_{P}$ denotes the Poincar\'e metric.
\end{proposition}

\begin{proof}
Let us prove it in the non-classical case. If such a constant $c$ does not
exist then there exists a sequence $\{f_{n}\}$ of orbifold morphism from $D$ to $(X/\Delta)$ such that
\begin{equation*}
\|f_{n}^{\prime }(0)\|_{\omega}>n.
\end{equation*}

Since $X$ is compact we may assume that $\{f_{n}(0)\}$
converges to a point $p\in X.$

Let $U$ be a complete hyperbolic neighborhood of $p$ in $X.$
Assume that there exists a positive number $r<1$ such that
$f_{n}(\Delta _{r})\subset U$ for $n\geq n_{0}.$ Then
$\{f_{n\left| \Delta _{r}\right. }:\Delta _{r}\rightarrow U\}$
would be relatively compact and woud have a subsequence which
converges to a holomorphic function from $\Delta _{r}$ to $U,$
which contradicts $\|f_{n}^{\prime }(0)\|_{\omega}>n.$

This means that for each positive integer $k,$ there exist a point
$z_{k}\in \Delta $ and an integer $n_{k}$ such that $\left|
z_{k}\right| <\frac{1}{k}$
and $f_{n_{k}}(z_{k})\notin U.$ Let $p_{k}=f_{n_{k}}(0)$ and $q_{k}=$ $%
f_{n_{k}}(z_{k}).$ By taking a subsequence we may assume that
$\{q_{k}\}$ converges to a point $q$ not in $U.$ Therefore we have
\begin{equation*}
d_{(X/\Delta)}(p_{k},q_{k})\leq d_{P}(0,z_{k})\rightarrow 0\text{ for }%
k\rightarrow \infty ,
\end{equation*}

and this contradicts the fact that $(X/\Delta)$ is hyperbolically imbedded
in $X$.

Conversely, let $\delta $ be the distance
function on $X$ induced by $\frac{1}{c}\omega$. Then
\begin{equation*} 
\delta \leq d_{(X/\Delta)}.
\end{equation*}

which implies obviously that $(X/\Delta)$ is hyperbolically imbedded in
$X$.
\end{proof}
\section{Algebraic hyperbolicity of orbifolds}
\subsection{The compact and the logarithmic setting}
In \cite{De95}, J.-P. Demailly introduced the concept of algebraic hyperbolicity for compact complex manifolds 
\begin{definition}
Let $X$ be a compact complex manifold and $\omega$ a hermitian metric on $X$. $X$ is algebraically hyperbolic if there exists $\varepsilon>0$ such that every compact irreducible curve $C \subset X$ satisfies
\begin{displaymath}
2g(\widetilde{C})-2\geq \varepsilon  \deg_{\omega}(C)
\end{displaymath}
where $g(\widetilde{C})$ is the genus of the normalization $\widetilde{C}$ of $C$ and $deg_{\omega}(C)=\int_{C}\omega$.
\end{definition}

Later, algebraic hyperbolicity was defined in the logarithmic setting by X. Chen in \cite{ch01} as follows
\begin{definition}
{\em Let $(X,D)$ be a log-manifold. For each reduced curve
$C\subset X$ that
meets $D$ properly, let $\nu :%
\widetilde{C}\rightarrow C$ be the normalization of $C.$ Then
$i(C,D)$ is the number of distinct points in the set $\nu
^{-1}(D)\subset \widetilde{C}.$}
\end{definition}

\begin{definition}\label{alghyp}
{\em A logarithmic variety $(X,D)$ is algebraically hyperbolic if
there exists a positive number $\varepsilon $ such that
\begin{equation*}
2g(\widetilde{C})-2+i(C,D)\geq \varepsilon \deg _{\omega }(C)
\end{equation*}
for all reduced and irreducible curves $C\subset X$ meeting $D$
properly where
$\widetilde{C}$ is the normalization of $C,$ $g(\widetilde{C})$ its genus and $%
\deg _{\omega }(C)=\int_{C}\omega$ with $\omega$ a hermitian
metric on $X$.}
\end{definition}

\subsection{The orbifold case}
We give the following definition which contains the two previous ones
\begin{definition}
Let $(X/\Delta)$ be a compact orbifold, $\omega$ a hermitian metric on $X$. $(X/\Delta)$ is algebraically (resp. classically algebraically) hyperbolic if there exists $\varepsilon>0$ such that for any non-constant orbifold (resp. classical orbifold) morphism $f: (C/\Delta') \rightarrow (X/\Delta)$, where $(C/\Delta')$ is an orbifold curve, 
\begin{displaymath}
\deg(K_{(C/\Delta')}):=\deg(K_{C}+\Delta')\geq \varepsilon \int_{C}f^*\omega.
\end{displaymath}
\end{definition}

The interest of such a notion is to provide a necessary condition for analytic hyperbolicity which is more tractable since it involves only algebraic curves instead of transcendental ones. So the first thing to do is to prove that analytic hyperbolicity implies algebraic hyperbolicity. In the compact setting this was done in \cite{De95} and for the logarithmic case in \cite{PR}. Here we will prove it in the more general setting of orbifolds.

\begin{theorem}
Let  $(X/\Delta)$ be an hyperbolic (resp. classically hyperbolic) compact orbifold hyperbolically imbedded in $X$. Then $(X/\Delta)$ is algebraically (resp. classically algebraically) hyperbolic.
\end{theorem}

\begin{proof}
Let $f: (C/\Delta') \rightarrow (X/\Delta)$ be a non-constant orbifold morphism where $(C/\Delta')$ is an orbifold curve. Since $(X/\Delta)$ is hyperbolic, $(C/\Delta')$ is hyperbolic. Therefore from theorem \ref{unif} we obtain that there is an orbifold morphism $\pi: D \rightarrow (C/\Delta')$. As explained above this gives us a singular metric $h_{(C/\Delta')}$ on $C\setminus \Delta'_{1}$ as the push-forward of the Poincar\'e metric by $\pi$. Moreover, since $(X/\Delta)$ is hyperbolically imbedded in $X$, we have $(f\circ \pi)^*\omega \leq c\ h_{P}$ with $c$ a positive constant. So, $f^*\omega \leq c\ h_{(C/\Delta')}$. Integrating and using theorem \ref{area}, we obtain $\frac{1}{2\pi c}\int_{C}f^*\omega \leq \deg(K_{(C/\Delta')}).$
\end{proof}

The following result generalizes a result of J.-P. Demailly \cite{De95} stating that any morphism $f:A \rightarrow X$ from an abelian variety to a compact algebraically hyperbolic manifold is constant.

\begin{theorem}
Let $(X/\Delta)$ be an algebraically (resp. classically algebraically) hyperbolic compact orbifold and $(\overline{A}/S)$ the logarithmic manifold associated to a semi-abelian variety $A$. Then any orbifold (resp. classically orbifold) morphism $f: (\overline{A}/S) \rightarrow (X/\Delta)$ is constant.
\end{theorem}

\begin{proof}
Let $C \subset A$ be a smooth curve with smooth compactification $\overline{C} \subset \overline{A}$. Let $m$ be a postive integer and consider the morphism $m_{A}: A \rightarrow A$ which is the multiplication by $m$ in $A$ and extends to a morphism $m_{\overline{A}}: \overline{A} \rightarrow \overline{A}$. We have the composition
\begin{displaymath}
f_{m}:\overline{C} \subset \overline{A} \overset{m_{\overline{A}}} \longrightarrow \overline{A} \overset{f}\longrightarrow X.
\end{displaymath}
Therefore we have an orbifold morphism
\begin{displaymath}
f_{m}:(\overline{C}/\Delta') \rightarrow (X/\Delta)
\end{displaymath}
where $\Delta'$ is the divisor induced on $\overline{C}$ by $S$. Let $L$ be an ample line bundle on $X$. We have
\begin{displaymath}
\deg(K_{\overline{C}}+\Delta')\geq \varepsilon \overline{C}.f_{m}^*L=\varepsilon m^2\overline{C}.f^*L.
\end{displaymath}
As $m$ can be as large as we want, this shows that $\overline{C}.f^*L=0$, so $f$ is constant on all curves in $A$ and therefore constant on $\overline{A}$.
\end{proof}

An interesting example of algebraically hyperbolic orbifold is given by the following theorem which generalizes results of \cite{ch01} and \cite{PR}
\begin{theorem}\label{algdeg}
Let $H_{1},\dots,H_{q}$ be very generic hypersurfaces of degrees $d_{i}$ in $\mathbb{P}^n$. Let $$\Delta = \sum_{1 \leq i \leq q} (1-\frac{1}{m_{i}})H_{i}.$$ Then for any non-constant orbifold morphism $f: (C/\Delta') \rightarrow (X/\Delta)$, where $(C/\Delta')$ is an orbifold curve, 
\begin{displaymath}
\deg(K_{(C/\Delta')}) \geq (\deg(\Delta)-2n)\deg(C).
\end{displaymath}
In particular if $\deg(\Delta)>2n$ then $(\mathbb{P}^n/\Delta)$ is algebraically hyperbolic.
\end{theorem}
\begin{proof}
Let $D=\sum_{i=1}^{q} H_{i}$, $C \subset \mathbb{P}^n$ be a reduced irreducible curve not contained in $D$ and $f:\tilde{C} \rightarrow C$ its normalization. Then it is sufficient to prove that $$\deg(K_{(\tilde{C}/\tilde{\Delta})}) \geq (\deg(\Delta)-2n)\deg(C)$$ for any orbifold $(\tilde{C}/\tilde{\Delta})$ such that $f:(\tilde{C}/\tilde{\Delta}) \rightarrow (\mathbb{P}^n/\Delta)$ is an orbifold morphism. First let us recall that in the logarithmic setting (i.e $m_{i}=\infty$ for all $i$) the result is known (see \cite{ch01} and \cite{PR}) i.e 
\begin{equation}
\label{ch}
2g(\widetilde{C})-2+i(C,D)\geq (d-2n) \deg(C)
\end{equation}
where $d=\sum_{i=1}^{q}d_{i}$.
There is a minimal orbifold structure on $\tilde{C}$ which makes $f$ an orbifold morphism and it is of course sufficient to prove the result for this orbifold. Let 
\begin{eqnarray*}
f^*(H_{j})&=&\sum_{i=1}^{i(C,D)}t_{i,j}p_{i},\\
f^*(D) & = & \sum_{i=1}^{i(C,D)}t_{i}p_{i},\\ 
\end{eqnarray*}
where the $p_{i}$ are the distinct points of $f^{-1}(D)$. Then if $\tilde{\Delta}=\sum_{i=1}^{i(C,D)}m_{i}'p_{i}$, the conditions for $f$ to be an orbifold morphism are $m_{i}'t_{i}\geq m_{j}$ for all $j\in \varphi(i)=\{1\leq k \leq i(C,D) / f(p_{i}) \in H_{k}\}$. Therefore the minimal orbifold structure is given by
$$\tilde{m_{i}}=\sup_{j\in \varphi(i)} \left\lceil\frac{m_{j}}{t_{i}}\right\rceil$$
where $\lceil k \rceil$ denotes the round up of $k$. So we have to prove
$$2g(\widetilde{C})-2+\sum_{i=1}^{i(C,D)}\left(1-\frac{1}{\tilde{m_{i}}}\right) \geq  (\deg(\Delta)-2n)\deg(C).$$
We have 
$$\sum_{i=1}^{i(C,D)}\left(1-\frac{1}{\tilde{m_{i}}}\right) \geq \sum_{i=1}^{i(C,D)}\left(1-\frac{t_{i}}{\sup_{j\in\varphi(i)}m_{j}}\right)\geq i(C,D)-\sum_{i=1}^{i(C,D)}\sum_{j=1}^{q}\frac{t_{i,j}}{m_{j}}.$$
Moreover
$$\sum_{i=1}^{i(C,D)}t_{i,j}=\deg(f^*H_{j})=C.H_{j}=\deg(C)d_{j}.$$
Therefore
$$\sum_{i=1}^{i(C,D)}\sum_{j=1}^{q}\frac{t_{i,j}}{m_{j}} = \deg(C)\left(\sum_{j=1}^{q}\frac{d_{j}}{m_{j}}\right).$$
From this inequality and \ref{ch} we obtain
$$2g(\widetilde{C})-2+\sum_{i=1}^{i(C,D)}\left(1-\frac{1}{\tilde{m_{i}}}\right) \geq (d-2n)\deg(C)- \deg(C)\left(\sum_{j=1}^{q}\frac{d_{j}}{m_{j}}\right).$$
And finally
$$2g(\widetilde{C})-2+\sum_{i=1}^{i(C,D)}\left(1-\frac{1}{\tilde{m_{i}}}\right) \geq (\deg(\Delta)-2n)\deg(C).$$
\end{proof}

This example is a motivation for the conjecture introduced in the next section.
As a corollary we obtain
\begin{corollary}\label{algdeg2}
Let $H_{1},\dots,H_{q}$ be very generic hypersurfaces of degrees $d_{i}$ in $\mathbb{P}^n$. Let $$\Delta = \sum_{1 \leq i \leq q} (1-\frac{1}{m_{i}})H_{i}$$ such that $\deg(\Delta)>2n$. Then every orbifold morphism $f:\bC \rightarrow(\mathbb{P}^n/\Delta)$ whose image is contained in an algebraic curve is constant.
\end{corollary}
\begin{proof}
Let $f(\bC) \subset C$. Then $f$ induces an orbifold morphism $f:\bC \rightarrow (C'/\Delta')$ where $C'$ is the normalization of $C$ and $\Delta'$ is the minimal orbifold structure making $(C'/\Delta') \rightarrow (\mathbb{P}^n/\Delta)$ an orbifold morphism. But by the previous theorem we have in particular that $(C/\Delta')$ is hyperbolic. This implies that $f$ is constant.
\end{proof}

\section{An orbifold Kobayashi's conjecture}
It is well known (see \cite{Ko98}) that the complement of $2n+1$ or more hyperplanes in general position in $\mathbb{P}^n$ is hyperbolic. In the orbifold setting we have 

\begin{theorem}\label{bhyp}
Let $H_{1}, H_{2},\dots,H_{q}$ be $q$ hyperplanes in general position in $\mathbb{P}^n$ with $q >2n$. Let $\Delta = \sum_{1 \leq i \leq q} (1-\frac{1}{m_{i}})H_{i}$ with $$\deg(\Delta)>q-\frac{q}{n}+1+\frac{1}{n}.$$ Then $(\mathbb{P}^n/\Delta)$ is hyperbolic and hyperbolically imbedded in $\mathbb{P}^n$.
\end{theorem}

We will prove this result using Nevanlinna theory so let us recall the usual notations. Let $E=\sum_{i=1}^{\infty} \nu_{i}z_{i}$ be a divisor on $\bC$ with distinct $z_{i} \in \bC$. Then we define the counting functions of $E$ truncated to $l\leq \infty$ by
\begin{eqnarray}
n_{l}(t,E) & = & \sum_{\{|z_{i}| < t\}} min{\{\nu_{i},l\}} \\
N^l(r,E) & = &  \int_{1}^r \frac{n_{l}(t,E)}{t}dt
\end{eqnarray}
If $D$ is a divisor on a complex space $X$ and $f:\bC \rightarrow X$ is a holomorphic map then
\begin{displaymath}
N_{f}^l(r,D)=N^l(r,f^*D).
\end{displaymath}
We denote $N_{f}(r,D):=N_{f}^{\infty}(r,D)$.

If $\omega$ is a $(1,1)$-form on $X$ then the order function with respect to $\omega$ is defined by
\begin{displaymath}
T_{f}(r,\omega)=  \int_{1}^r \Big(\int_{|z|<t} f^*\omega\Big) \frac{dt}{t}.
\end{displaymath}
The defect is defined by
\begin{displaymath}
\delta^l(f,D)=\underset{r \rightarrow \infty} \liminf \Big(1-\frac{N_{f}^l(r,D)}{T_{f}(r,c_{1}(D))}\Big).
\end{displaymath}

To state the result we need, let us recall the notion of hyperplanes in subgeneral position. Let $N\geq n$ and $q\geq N+1$. We consider $q$ hyperplanes $H_{1}, H_{2},\dots,H_{q}$ in $\mathbb{P}^n$, which are given by $$H_{j}: <Z,A_{j}>=0$$ for non-zero vectors $A_{j}$ in $\bC^{n+1}$. Then we say that  $H_{1}, H_{2},\dots,H_{q}$ are in $N$-subgeneral position if for every $1 \leq i_{0} \leq \dots \leq i_{N} \leq q$  $$span(A_{i_{0}},\dots,A_{i_{N}})=\bC^{n+1}.$$
We have the following generalized defect relation of Cartan due to Nochka \cite{No83}:
\begin{theorem}\label{No}
Let $H_{1}, H_{2},\dots,H_{q}$ be $q$ hyperplanes in $N$-subgeneral position in $\mathbb{P}^n$. Then for any holomorphic map $f:\bC \rightarrow \mathbb{P}^n$ that is non-linearly degenerate, we have
\begin{displaymath}
\sum_{i=1}^{q}\delta^n(f,H_{i})\leq 2N-n+1.
\end{displaymath}
\end{theorem}

As a corollary we obtain
\begin{theorem}\label{deg}
Let $H_{1}, H_{2},\dots,H_{q}$ be $q$ hyperplanes in general position in $\mathbb{P}^n$ with $q >2n$. Let $\Delta = \sum_{1 \leq i \leq q} (1-\frac{1}{m_{i}})H_{i}$ with $$\deg(\Delta)>q-\frac{q}{n}+1+\frac{1}{n}.$$ Then every orbifold morphism $f:\bC \rightarrow (\mathbb{P}^n/\Delta)$ is constant.
\end{theorem}
\begin{proof}
Suppose $\bP^l \subset \bP^n$ contains $f(\bC)$. The intersections of $\bP^l$ with the $H_{i}$ are in $n$-subgeneral position in $\bP^l$.
By the First Main Theorem of Nevanlinna theory we have 
\begin{displaymath}
T_{f}(r,c_{1}(H_{i}))\geq N_{f}(r,H_{i})+C
\end{displaymath}
where $C$ is a constant.
Since $f^*H_{i}$ has multiplicity at least $m_{i}$ at every point of $f^{-1}H_{i}$ we have
\begin{displaymath}
N_{f}(r,H_{i})\geq \frac{m_{i}}{l}N^l_{f}(r,H_{i}).
\end{displaymath}
Therefore
\begin{displaymath}
\delta^l(f,H_{i})\geq 1-\frac{l}{m_{i}}.
\end{displaymath}
Now
\begin{displaymath}
\sum_{i=1}^{q}\delta^l(f,H_{i}) \geq \sum_{i=1}^{q} (1-\frac{l}{m_{i}})=l\deg(\Delta)-(l-1)q.
\end{displaymath}
Therefore
\begin{displaymath}
\sum_{i=1}^{q}\delta^l(f,H_{i}) > q(1-\frac{l}{n})+l+\frac{l}{n}.
\end{displaymath}
From theorem \ref{No}, we deduce that the conclusion holds provided
\begin{displaymath}
q(1-\frac{l}{n})+l+\frac{l}{n} \geq 2n-l+1
\end{displaymath}
for $l=1,\dots,n$. This condition is clearly satisfied since
\begin{displaymath}
q(1-\frac{l}{n})\geq (2n+1)(1-\frac{l}{n})=2n+1-2l-\frac{l}{n}.
\end{displaymath}
\end{proof}

Now, we can prove theorem \ref{bhyp}.
\begin{proof}
Suppose $(\bP^n/\Delta)$ is not hyperbolically imbedded in $\mathbb{P}^n$. We shall need a slight refinement of theorem \ref{brody2}. There is a sequence of orbifold morphisms $f_{n}:D \rightarrow (\bP^n/\Delta)$ such that $lim ||f'_{n}(0)||=+\infty$. Thanks to Brody reparametrization, we obtain a sequence of orbifold morphisms $g_{n}:D(0,r_{n})\rightarrow (\bP^n/\Delta)$, with $r_{n} \rightarrow +\infty$, converging to a holomorphic map $f: \bC \rightarrow \bP^n$ which is either a non constant orbifold morphism $f: \bC \rightarrow (\bP^n/\Delta)$ or a non-constant holomorphic map $f: \bC \rightarrow supp(\Delta)$. 

The first case is not possible thanks to theorem \ref{deg}.

Consider a partition of indices $\{1,2,\dots,q\}=I\cup J$, and let $L_{I}=\cap_{i \in I} H_{i}$. If $I$ contains $k$ elements, $L_{I}$ is an $n-k$ dimensional linear subspace. The intersections $Z_{j}=H_{j}\cap L_{I}$ are $q-k$ hyperplanes in general position in $L_{I}$. 

The observation here is that in the second case there is a partition $\{1,2,\dots,q\}=I\cup J$ such that $f$ is an orbifold morphism from $\bC$ to $(L_{I}/\Delta')$ where $\Delta'=\sum_{i\in J} (1-\frac{1}{m_{i}})Z_{i}$. Indeed, the sequence $g_{n}:D(0,r_{n})\rightarrow (\bP^n/\Delta)$ can be seen as a sequence of orbifold morphisms $g_{n}:D(0,r_{n})\rightarrow (\bP^n/\Delta_{J})$, where $\Delta_{J}=\sum_{i\in J} (1-\frac{1}{m_{i}})H_{i}$ since $\Delta_{J} \leq \Delta$. Therefore (see \cite{CW}) it converges to a map $f$ which is either an orbifold morphism  from $\bC$ to $(\bP^n/\Delta')$ or verifies $f(\bC) \subset supp(\Delta')$.

Using again theorem \ref{deg}, we see that the conclusion holds provided
\begin{eqnarray*}
q-k & > & 2(n-k),\\
\deg(\Delta_{J}) & > & q-k-\frac{q-k}{n-k}+1+\frac{1}{n-k} 
\end{eqnarray*}
for $k=1,\dots,n-1$.
The first condition is satisfied since $q>2n$.
For the second one, suppose it is not true then
$$\deg(\Delta)=\deg(\Delta_{I})+\deg(\Delta_{J})\leq  q-k-\frac{q-k}{n-k}+1+\frac{1}{n-k} +k.$$
Therefore
$$\deg(\Delta)-(q-\frac{q}{n}+1+\frac{1}{n}) \leq \frac{q-1}{n}-\frac{q-k-1}{n-k}=\frac{k(n+1-q)}{n(n-k)}<0.$$
So we obtain a contradiction.
\end{proof}

These results and theorem \ref{algdeg} suggest the following generalization of a conjecture of S. Kobayashi (see \cite{Ko98}):

\begin{conjecture}\label{KO}
Let $H_{0}, H_{1},\dots,H_{q}$ be generic hypersurfaces in $\mathbb{P}^n$. Let $$\Delta = \sum_{0 \leq i \leq q} (1-\frac{1}{m_{i}})H_{i}$$ with $\deg(\Delta)>2n$. Then $(\mathbb{P}^n/\Delta)$ is hyperbolic and hyperbolically imbedded in $\mathbb{P}^n$.
\end{conjecture}
We will prove some more results towards this conjecture in the next section.

\section{Orbifold jet differentials and applications}
\subsection{Symmetric differentials}
Let $(X/\Delta)$ be a smooth orbifold i.e $\lceil \Delta \rceil := supp(\Delta)$ is a normal crossing divisor. Let $(x_{1},\dots,x_{n})$ be local coordinates such that $\Delta$ has equation $$x_{1}^{(1-\frac{1}{m_{1}})}\dots x_{n}^{(1-\frac{1}{m_{n}})}=0.$$
Following \cite{C07} we can define sheaves of differential forms on orbifolds
\begin{definition}
For $N$ a positive integer, $S^{N}\Omega_{(X/\Delta)}$ is the locally free subsheaf of $S^{N}\Omega_{X}(log\lceil \Delta \rceil)$ generated by the elements $$x_{1}^{\lceil{\frac{\alpha_{1}}{m_{1}}}\rceil}\dots x_{n}^{\lceil{\frac{\alpha_{n}}{m_{n}}}\rceil}\left(\frac{dx_{1}}{x_{1}}\right)^{\alpha_{1}}\dots\left(\frac{dx_{n}}{x_{n}}\right)^{\alpha_{n}}$$
such that $\sum \alpha_{i}=N$, where $\lceil k \rceil$ denotes the round up of $k$.
\end{definition}

As an immediate consequence of the definition we have
\begin{proposition}\cite{C07}
Let $(X/\Delta)$ and $(Y/\Delta')$ be smooth orbifolds and $f: (X/\Delta) \rightarrow (Y/\Delta')$ an orbifold morphism. Then $$f^*(S^{N}\Omega_{(Y/\Delta')}) \subset S^{N}\Omega_{(X/\Delta)}.$$
\end{proposition}

\subsection{Applications to hyperbolicity}
It is well known in the compact and logarithmic cases that the existence of global symmetric differential forms vanishing on an ample divisor provide differential equations for entire curves. This can be generalized to the orbifold case as shown in \cite{CP}
\begin{theorem}(Campana-Paun)\label{vt}
Let  $(X/\Delta)$ be a smooth compact orbifold, $A$ an ample line bundle on $X$ and $\omega \in H^0(X,S^{N}\Omega_{(X/\Delta)}\otimes A^{-1})$. Then for any orbifold morphism $f:\bC \rightarrow (X/\Delta)$
\begin{displaymath}
\omega(f')\equiv 0.
\end{displaymath}
\end{theorem}
Here we would like to illustrate these ideas initiated in \cite{CP}. We provide some results towards conjecture \ref{KO} which  generalize \cite{EG} where the logarithmic case was studied. 

\begin{theorem}\label{gs}
Let $(X/\Delta)$ be a smooth projective orbifold surface of general type, where $\Delta=\sum_{i=1}^{n}(1-\frac{1}{m_{i}})C_{i}$ is the decomposition into irreducible components, and $A$ an ample line bundle on $X$. Suppose that $g_{i}:=g(C_{i}) \geq 2$, $h^0(C_{i},\mathcal{O}_{C_{i}}(C_{i})) \neq 0$ for all $i$ and
that the logarithmic Chern classes of $(X, \lceil \Delta \rceil)$ verify
\begin{equation}
\label{eq}
\overline{c_{1}}^2-\overline{c_{2}}-\sum_{i=1}^{n}\frac{1}{m_{i}}(2g_{i}-2+\sum_{j \neq i}C_{i}C_{j})>0
\end{equation}
then $$H^0(X,S^{N}\Omega_{(X/\Delta)}\otimes A^{-1})\neq 0$$
for $N$ large enough.
\end{theorem}

\begin{proof}
Let $F:=S^{N}\Omega_{X}(\log\lceil \Delta \rceil)/S^{N}\Omega_{(X/\Delta)}$ be the quotient sheaf. $F$ is supported on $C=\bigcup C_{i}$ and 
\begin{eqnarray*}
h^0(X,S^{N}\Omega_{(X/\Delta)}) & \geq& h^0(X,S^{N}\Omega_{X}(\log\lceil \Delta \rceil))-h^0(C,F) \\
& \geq &  h^0(X,S^{N}\Omega_{X}(log\lceil \Delta \rceil))- \sum_{i=1}^{n}h^0(C_{i},F_{|C_{i}}).
\end{eqnarray*}
$F_{|C_{i}}$ has a natural filtration
$$ F_{N} \subset F_{N-1} \subset \dots \subset F_{1} \subset F_{0}=F_{|C_{1}}$$ 
such that $$F_{j}/F_{j+1}=S^{N-j}\Omega_{(C_{i}/\Delta_{i})}\otimes G_{j}$$
where $\Delta_{i}:=\sum_{j \neq i}(1-\frac{1}{m_{j}})\mathcal{O}_{C_{i}}(C_{j})$ is the divisor induced by $\Delta$ on $C_{i}$ and $G_{j}$ is a locally free sheaf of rank $r_{j}:= \lceil \frac{j}{m_{i}} \rceil$ admitting a filtration
$$ H_{1}\subset H_{1} \subset \dots \subset H_{r_{j}}=G_{j}$$
where $$H_{i}/H_{i-1}=[N_{C_{i}}^*]^{\otimes i}$$
with $N_{C_{i}}$ denoting the normal bundle of $C_{i}$.
We have 
$$S^{N-j}\Omega_{(C_{i}/\Delta_{i})}=\lfloor (N-j)(K_{C_{i}}+\Delta_{i}) \rfloor.$$
Since $[N_{C_{i}}]=\mathcal{O}_{C_{i}}(C_{i})$ and $h^0(C_{i},\mathcal{O}_{C_{i}}(C_{i})) \neq 0$ we obtain
\begin{eqnarray*}
h^0(C_{i},F_{|C_{i}})& \leq & \sum_{j=0}^{N}{h^0(C_{i},F_{j}/F_{j+1})}\\
&\leq &  h^0(C_{i},N(K_{C_{i}}+\sum_{j \neq i}\mathcal{O}_{C_{i}}(C_{j}))+\sum_{j=1}^{N}{r_{j}h^0(C_{i},(N-j)(K_{C_{i}}+\sum_{j \neq i}\mathcal{O}_{C_{i}}(C_{j}))}.
\end{eqnarray*}
Since  $g_{i} \geq 2$ we obtain for $N-j \geq 1$
$$h^0(C_{i},(N-j)(K_{C_{i}}+\sum_{j \neq i}\mathcal{O}_{C_{i}}(C_{j})) \leq (N-j)(2g_{i}-2+\sum_{j \neq i}C_{i}C_{j})-g_{i}+1.$$
Now we suppose $N=qm_{i}$ is a multiple of $m_{i}$. We write $j=(h-1)m_{i}+k$ with $1 \leq k \leq m_{i}$ so $r_{j}:=h$. Therefore we obtain
\begin{eqnarray*}
h^0(C_{i},F_{|C_{i}})& \leq & m_{i}q(2g_{i}-2+\sum_{j \neq i}C_{i}C_{j})-g_{i}+1+\sum_{h=1}^{q}{m_{i}(hm_{i}(q-h+1)(2g_{i}-2+\sum_{j \neq i}C_{i}C_{j})-g_{i}+1)}\\
& \leq & \frac{1}{6}m_{i}^2q^3(2g_{i}-2+\sum_{j \neq i}C_{i}C_{j})+O(q^2).
\end{eqnarray*}
Therefore we obtain if $N=qm_{1}\dots m_{n}$
\begin{displaymath}
\sum_{i=1}^{n}h^0(C_{i},F_{|C_{i}})\leq \frac{1}{6}(m_{1}\dots m_{n}q)^3\big(\sum_{i=1}^{n}\frac{1}{m_{i}}(2g_{i}-2+\sum_{j \neq i}C_{i}C_{j})\big)+O(q^2).
\end{displaymath}
On the other hand
\begin{displaymath}
\chi(X,S^{N}\Omega_{X}(\log\lceil \Delta \rceil)=\frac{(m_{1}\dots m_{n}q)^3}{6}(\overline{c_{1}}^2-\overline{c_{2}})+O(q^2).
\end{displaymath}
Using Bogomolov's vanishing theorem \cite{Bog} for the $h^2$ term we obtain
\begin{displaymath}
h^0(X,S^{N}\Omega_{X}(\log\lceil \Delta \rceil)) \geq \frac{(m_{1}\dots m_{n}q)^3}{6}(\overline{c_{1}}^2-\overline{c_{2}})+O(q^2).
\end{displaymath}
So finally we have
\begin{eqnarray*}
h^0(X,S^{N}\Omega_{(X/\Delta)}) & \geq & \frac{(m_{1}\dots m_{n}q)^3}{6}\big(\overline{c_{1}}^2-\overline{c_{2}}-\sum_{i=1}^{n}\frac{1}{m_{i}}(2g_{i}-2+\sum_{j \neq i}C_{i}C_{j})\big)+O(q^2).
\end{eqnarray*}
\end{proof}

As a first application, we can generalize a result due to Bogomolov \cite{Bog77} about the finiteness of rational and elliptic curves on surfaces of general type with $c_{1}^2>c_{2}$.  
\begin{theorem}\label{Bog}
Let $(X/\Delta)$ be a smooth compact orbifold surface of general type with the same hypotheses as in theorem \ref{gs}. Then there are only finitely many special curves $C\subset (X/\Delta)$ i.e images of non-constant orbifold morphism $\nu: (C'/\Delta') \rightarrow (X/\Delta)$ where $(C'/\Delta')$ is an orbifold curve with $\deg(K_{(C'/\Delta')}\leq 0$.
\end{theorem}
\begin{proof}
Let $Y=\bP(T_{X}(-\log\lceil \Delta \rceil))$ be the projectivization of the logarithmic tangent bundle. From theorem \ref{gs} there is a global section $\omega \in H^0(X,S^{N}\Omega_{(X/\Delta)}\otimes A^{-1})$ which can be seen as a holomorphic section of $\mathcal{O}_{Y}(N)\otimes\pi^*A^{-1}$, where $\pi : Y \rightarrow X$ denotes the canonical projection. Now, from theorem \ref{vt}, the lifts of the special curves must lie in an irreducible component $Z \subset Y$ of the zeros of $\omega$. Let $V \subset T_{Y}(-\log \pi^*\lceil \Delta \rceil)$ be the subbundle defined by
\begin{displaymath}
V_{x,[v]}=\{\xi \in T_{Y}(-\log \pi^*\lceil \Delta \rceil); (\pi)_{*}\xi\in\bC.v\}.
\end{displaymath}
Then $V$ defines on the desingularization $\tilde{Z}$ of $Z$ an algebraic foliation by curves, such that the tangent bundle to the leaves is given by $T_{Z}\cap V$. The lifts of the special curves are leaves of this foliation. Now, a theorem of Jouanolou \cite{Jou78} implies that if there are an infinite number of such curves then there is a meromorphic fibration $\tilde{Z} \rightarrow S$ from $\tilde{Z}$ to a curve such that the leaves correspond to the fibers of the fibration. Let $\tilde{\Delta} \subset \tilde{Z}$ be the divisor above $\Delta$. Then $(\tilde{Z}/\tilde{\Delta})$ is of general type. But we cannot have a fibration $(\tilde{Z}/\tilde{\Delta}) \rightarrow S$ with special generic fiber $(\tilde{Z}/\tilde{\Delta})_{s}$ and at the same time $(\tilde{Z}/\tilde{\Delta})$ of general type (see corollary 7.14 of \cite{C07}).
\end{proof}

Now, we turn to the transcendental case. We will use the following result which is a consequence of deep works of McQuillan on foliations of surfaces (see \cite{McQ0}, \cite{McQ}, \cite{McQ2} and also \cite{EG} for a weaker version). For the convenience of the reader we will give a proof, refering to the above mentioned articles for details.
\begin{theorem}\label{McQ}
Let $(X,D)$ be a smooth logarithmic projective surface of log general type and $f : \bC \rightarrow X$ a non algebraically degenerate entire curve. Let $f_{[1]}: \bC \rightarrow \bP(T_{X}(-\log D))$ denote the canonical lifting of $f$. Suppose $f_{[1]} (\bC)$ is contained in a divisor in $\bP(T_{X}(-\log D))$. Then 
$$T_{f}(r,c_{1}(K_{X}+D)) \leq N^{1}_{f}(r,D)+\epsilon T_{f}(r,c_{1}(H))||_{\epsilon}$$
for some ample line bundle $H$ on $X$ and the notation $||_{\epsilon}$ meaning that the inequality holds for any $\epsilon >0$, for $r$ outside a subset of finite measure depending on $\epsilon$.
\begin{proof}
Let $S \subset \bP(T_{X}(-\log D))$ denote the surface which contains $f(\bC)$ and $\pi: S \rightarrow X$ the canonical projection. As already explained in the proof of theorem \ref{Bog}, there is a canonical foliation $\mathcal{F}$ on $S$ such that $f_{[1]}: \bC \rightarrow \bP(T_{X}(-\log D))$ is a leaf of  $\mathcal{F}$. After some blow ups we obtain a foliated smooth surface $(\tilde{S}, \tilde{D}, \tilde{\mathcal{F}}) \rightarrow (S, \pi^{-1}(D), \mathcal{F})$, i.e $\tilde{S}$ is smooth,  $\tilde{D}$ is normal crossing and $\tilde{\mathcal{F}}$ has reduced singularities. Let $\tilde{D}=C+B$ where $C$ is the invariant part of $\tilde{D}$ by $\tilde{\mathcal{F}}$. We have an exact sequence
$$0 \rightarrow \mathcal{N}^{*}(C) \rightarrow T^*_{\tilde{S}}(\log \tilde{D}) \rightarrow K_{\tilde{\mathcal{F}}}(B).\mathcal{I}_{Z} \rightarrow 0,$$
where $\mathcal{I}_{Z}$ is an ideal supported on the singularity set $Z$ of $\mathcal{F}$. Now, we apply the logarithmic tautological inequality of McQuillan (see \cite{McQ01} and \cite{Vo}) which gives
$$T_{\tilde{f}_{[1]}}(r,c_{1}(L)) \leq N^{1}_{\tilde{f}}(r, \tilde{D}) + \epsilon T_{f}(r,c_{1}(H))||_{\epsilon},$$
where $L=\mathcal{O}_{\bP(T_{\tilde{S}}(-\log \tilde{D}))}(1)$, $\tilde{f}$ and $\tilde{f}_{[1]}$ are the lifts of $f$. Moreover, we have
$$L_{|Y}=p^*K_{\tilde{\mathcal{F}}}(B) \otimes\mathcal{O}(-E),$$
where $L_{|Y}$ denotes the restriction of $L$ to the graph $Y$ of the foliation, $p:Y \rightarrow \tilde{S}$ the projection and $E$ is the total exceptional divisor. 
Therefore we obtain
$$T_{\tilde{f}}(r,c_{1}(K_{\tilde{\mathcal{F}}}+B))\leq N^{1}_{\tilde{f}}(r, \tilde{D}) + \epsilon T_{f}(r,c_{1}(H))||_{\epsilon}.$$
Indeed, it is proved in \cite{McQ0} (see also \cite{McQ}, \cite{McQ2}, \cite{Br} and \cite{EG}), using sequence of blow ups on $\tilde{S}$, that $T_{\tilde{f}_{[1]}}(r,c_{1}(E))\leq \epsilon T_{f}(r,c_{1}(H))||_{\epsilon}.$ It is also proved (see  \cite{McQ0}, \cite{McQ}, \cite{McQ2}, \cite{Br} and \cite{EG}) that $T_{\tilde{f}}(r,c_{1}(\mathcal{N}^{*}(C))) \leq \epsilon T_{f}(r,c_{1}(H))||_{\epsilon}.$
Therefore we obtain
$$T_{\tilde{f}}(r,c_{1}(K_{\tilde{S}}+\tilde{D})) \leq N^{1}_{\tilde{f}}(r, \tilde{D}) + \epsilon T_{f}(r,c_{1}(H))||_{\epsilon},$$
and finally
$$T_{f}(r,c_{1}(K_{X}+D)) \leq N^{1}_{f}(r,D)+\epsilon T_{f}(r,c_{1}(H))||_{\epsilon}.$$
\end{proof}

\end{theorem}

Now, we can prove
\begin{theorem}\label{MQdeg}
Let $(X/\Delta)$ be a smooth compact orbifold surface of general type with the same hypotheses as in theorem \ref{gs}. Then every orbifold morphism $f: \bC \rightarrow (X/\Delta)$ is algebraically degenerate.
\end{theorem}
\begin{proof}
Suppose $f: \bC \rightarrow (X/\Delta)$ is non algebraically degenerate. Let $D:=\lceil \Delta \rceil$. By the hypotheses, $f_{[1]} (\bC)$ is contained in a divisor in $\bP(T_{X}(-\log D))$. Then we can apply theorem \ref{McQ} and we have
$$T_{f}(r,c_{1}(K_{X}+D)) \leq N^{1}_{f}(r,D)+\epsilon T_{f}(r,c_{1}(H))||_{\epsilon}.$$
Moreover $$T_{f}(r,c_{1}(K_{X}+\Delta)) =T_{f}(r,c_{1}(K_{X}))+\sum_{i=1}^{n}(1-\frac{1}{m_{i}})T_{f}(r,c_{1}(C_{i})).$$
We have $m_{i}N^{1}_{f}(r,C_{i})\leq N_{f}(r,C_{i})$ and by the First Main Theorem $N_{f}(r,C_{i}) \leq T_{f}(r,c_{1}(C_{i})) + O(1).$ Therefore 
$$T_{f}(r,c_{1}(C_{i}))-N^{1}_{f}(r,C_{i})\geq (1-\frac{1}{m_{i}})T_{f}(r,c_{1}(C_{i}))+O(1).$$
So we obtain
$$T_{f}(r,c_{1}(K_{X}+\Delta))\leq T_{f}(r,c_{1}(K_{X}))+\sum_{i=1}^{n}(T_{f}(r,c_{1}(C_{i}))-N^{1}_{f}(r,C_{i}))+O(1).$$
And
$$T_{f}(r,c_{1}(K_{X}+\Delta))\leq T_{f}(r,c_{1}(K_{X}+D))-N^{1}_{f}(r,D)+O(1),$$
which gives
$$T_{f}(r,c_{1}(K_{X}+\Delta))\leq \epsilon T_{f}(r,c_{1}(H))||_{\epsilon}.$$
This is a contradiction since  $(X/\Delta)$ is of general type, and $f$ is algebraically degenerate.
\end{proof}

\begin{remark}
These results generalize \cite{EG} where all multiplicities are infinite.
\end{remark}

As a consequence of theorems \ref{Bog} and \ref{MQdeg} we obtain the proof of 

\bigskip

\noindent {\bf Theorem \ref{theo1}.} \textit{
Let $(X/\Delta)$ be a smooth compact orbifold surface of general type with the same hypotheses as in theorem \ref{gs}. Then there exists a proper subvariety $Y\subsetneq X$ such that every orbifold morphism $f: \bC \rightarrow (X/\Delta)$ verifies $f(\bC) \subset Y$.}

\subsubsection{Degeneracy of holomorphic curves with ramification}

As an application we obtain the following theorem
\begin{theorem}
Let $C_{i}, 1\leq i \leq 2$, be two smooth curves in $X=\bP_{2}$ of degree $d_{i}\geq 4$ with normal crossings. Let $\Delta=(1-\frac{1}{m_{1}})C_{1}+(1-\frac{1}{m_{2}})C_{2}$, and $d=d_{1}+d_{2}$. If 
\begin{equation}
\label{eq1}
deg(\Delta)>\frac{d_{1}^2+d_{2}^2+d_{1}d_{2}-6}{d-3}
\end{equation}
then every orbifold morphism $f: \bC \rightarrow (X/\Delta)$ is algebraically degenerate.
Moreover, if the curves $C_{i}$ are very generic, then $(X/\Delta)$ is hyperbolic.
\end{theorem}

\begin{proof}
First we verify that condition \ref{eq} is satisfied. We compute everything in terms of the degrees $d_{1}\leq d_{2}$
\begin{eqnarray*}
\overline{c_{1}}^2-\overline{c_{2}}-\frac{1}{m_{1}}(2g_{1}-2+d_{1}d_{2})-\frac{1}{m_{2}}(2g_{2}-2+d_{1}d_{2})=deg(\Delta)(d-3)-(d_{1}^2+d_{2}^2+d_{1}d_{2}-6).
\end{eqnarray*}
So, if condition \ref{eq1} is satisfied, we can apply theorem \ref{theo1} and obtain the algebraic degeneracy of $f$. If $\overline{c_{1}}^2-\overline{c_{2}}>0$ then $d_{1}\geq 5$ or $d_{1}\geq 4$ and $d_{2} \geq 7$ (see \cite{Rou03}). Therefore $\deg(\Delta) >4$ and, if the curves are very generic, from corollary \ref{algdeg2} we obtain that $f$ is constant.
\end{proof}

\begin{example}
Let $C_{i}, 1\leq i \leq 2$ be two smooth curves in $\bP_{2}$ of degree $5$ with normal crossings. Let $\Delta=(1-\frac{1}{70})C_{1}+(1-\frac{1}{71})C_{2}$. Then every orbifold morphism $f: \bC \rightarrow (X/\Delta)$ is algebraically degenerate. If the curves $C_{i}$ are very generic, then $(\bP_{2}/\Delta)$ is hyperbolic.
\end{example}

\subsubsection{Weakly special manifolds with degenerate entire curves}
Let us recall that a complex projective manifold is said to be weakly-special if none of its finite etale covers has a dominant rational map to a positive-dimensional manifold of general type. Bogomolov and Tschinkel constructed in \cite{BT} examples of algebraic threefolds which are weakly-special but non-special. These examples are simply connected and come with an elliptic fibration $\varphi: X \rightarrow B$, where $B$ is a surface with $\kappa(B)=1$. This fibration is of general type because it has multiple fibres of multiplicity $m\geq 2$ over a smooth curve $D$ such that $\kappa(B,K_{B}+(1-\frac{1}{m})D)=2$. But these multiple fibers cannot be eliminated by an etale cover of $X$, which is simply connected. We remark that our results enable us to simplify the proof and generalize the results of \cite{CP} where it is proved that in some of the examples of  \cite{BT} all entire curves are degenerate. Indeed, using these notations, an immediate consequence of theorem \ref{theo1} is
\begin{theorem}
Let $X$ be a Bogomolov-Tshinkel example. If $$\overline{c_{1}}^2(B,D)-\overline{c_{2}}(B,D)-\frac{1}{m}(2g(D)-2)>0$$ then there exists $\Gamma \subsetneq B$ such that for any entire curve $h:\bC \rightarrow X$, $\phi \circ h: \bC \rightarrow B$ is either a point or contained in $\Gamma$. 
\end{theorem}
\begin{proof}
We just remark that $\phi \circ h: \bC \rightarrow (B,(1-\frac{1}{m})D)$ is an orbifold morphism and apply theorem \ref{theo1}.
\end{proof}

\begin{remark}
The use of McQuillan's results in the logarithmic setting simplify the arguments used in \cite{CP} to obtain the algebraic degeneracy. Indeed, the authors need additional technical hypotheses to use the compact version of McQuillan's results.
\end{remark}

\subsection{Higher order jet differentials}

More generally, we generalize jet differentials from the compact setting (see \cite{GG80} and \cite{De95}) and the logarithmic setting (see \cite{DL96}) to the orbifold setting. 

Recall that if $X$ is a compact complex manifold, in \cite{GG80} Green and Griffiths have introduced the vector bundle of jet differentials of order $k$ and degree $m$, $E_{k,m}^{GG}\Omega_{X}\rightarrow X$ whose fibers are complex valued polynomials $Q(f^{\prime },f^{\prime \prime },\dots,f^{(k)})$ on the fibers of $J_{k}X$ of weight $m$ for the action of  $\mathbb{C}^{\ast }$:
\begin{equation*}
Q(\lambda f^{\prime },\lambda ^{2}f^{\prime \prime },\dots,\lambda
^{k}f^{(k)})=\lambda ^{m}Q(f^{\prime },f^{\prime \prime },\dots,f^{(k)})
\end{equation*}
for all $\lambda \in \mathbb{C}^{\ast }$ and $(f^{\prime },f^{\prime \prime
},\dots,f^{(k)})\in J_{k}X.$

If $(X,D)$ is a logarithmic manifold, i.e $X$ is a compact complex manifold and $D=\sum_{i} D_{i}$ is a normal crossing divisor, the vector bundle of logarithmic jet differentials of order $k$ and degree $m$, $E_{k,m}^{GG}\Omega_{(X,D)}\rightarrow X$ consists of polynomials satisfying the same weight condition $Q(f^{\prime },f^{\prime \prime },\dots,f^{(k)})$ in the derivatives of $f$ and in the derivatives of $\log(s_{j}(f))$ on $D_{j}=\{s_{j}=0\}$ locally.

Let $(X/\Delta)$ be a smooth orbifold. Let $(x_{1},\dots,x_{n})$ be local coordinates such that $\Delta$ has equation $$x_{1}^{(1-\frac{1}{m_{1}})}\dots x_{n}^{(1-\frac{1}{m_{n}})}=0.$$
\begin{definition}
For $N$ a positive integer, $E_{k,N}\Omega_{(X/\Delta)}$ is the locally free subsheaf of $E_{k,N}^{GG}\Omega_{(X,\lceil \Delta \rceil)}$ generated by the elements
$$\prod_{1\leq i \leq n}x_{i}^{\lceil \frac{\alpha_{i,1}}{m_{i}}\rceil}\left(\frac{dx_{i}}{x_{i}}\right)^{\alpha_{i,1}}\dots\prod_{1 \leq i \leq n}x_{i}^{\lceil \frac{k\alpha_{i,k}}{m_{i}}\rceil}\left(\frac{d^kx_{i}}{x_{i}}\right)^{\alpha_{i,k}}$$
such that $|\alpha_{1}|+2|\alpha_{2}|+\dots+k|\alpha_{k}|=N$ where $|\alpha_{i}|=\sum_{j}\alpha_{j,i}$.
\end{definition}

Now theorem \ref{vt} generalizes as
\begin{theorem}
Let  $(X/\Delta)$ be a smooth compact orbifold, $A$ an ample line bundle on $X$ and $P \in H^0(X,E_{k,N}\Omega_{(X/\Delta)}\otimes A^{-1})$. Then for any orbifold morphism $f:\bC \rightarrow (X/\Delta)$
\begin{displaymath}
P(f)\equiv 0.
\end{displaymath}
\end{theorem}

\begin{proof}
The proof goes along the same lines as in the classical setting using the logarithmic derivative lemma (see \cite{Siu}, \cite{W99}, \cite{CP}) which we summarize for the convenience of the reader.
$P(f)$ is a holomorphic section of $f^*A^{-1}$. Suppose it does not vanish identically. Let $\omega=\Theta_{h}(A)$, then by the Poincar\'e-Lelong equation
$$i\partial \overline\partial \log ||P(f)||^2_{h^{-1}} \geq f^*\omega.$$ Therefore
$$T_{f}(r,\omega)\leq \int_{1}^{r}\frac{dt}{t} \int_{|z|<t} i\partial \overline\partial \log ||P(f)||^2_{h^{-1}}$$
and from Jensen formula
$$\int_{0}^{2\pi}\log^{+} ||P(f)||_{h^{-1}}d\theta \geq T_{f}(r,\omega)+\mathcal{O}(1).$$
Finally the logarithmic derivative lemma gives
$$\int_{0}^{2\pi}\log^{+} ||P(f)||_{h^{-1}}d\theta \leq \mathcal{O}(\log(r)+\log(T_{f}(r,\omega))$$
outside a set of finite Lebesgue measure in $[0,+\infty[$.
This gives a contradiction.
\end{proof}

\begin{remark}
In a sequel of this paper we will give some applications of these higher order jet differentials and develop an approach using stacks.
\end{remark}

\section{Measure hyperbolicity and orbifolds of general type}
\subsection{Kobayashi-Ochiai's extension theorems}
Let us recall some results about holomorphic mappings into orbifolds of general type.

In \cite{KO75} it was established that a meromorphic map of maximal rank from a dense Zariski open subset $U$ of a complex manifold $V$ to a compact variety of general type extends meromorphically to $V$.

The main results of \cite{Sa74} can be formulated in the orbifold setting as follows
\begin{theorem}(Sakai)
Let $(X/\Delta)$ be a smooth orbifold of general type where $X$ is a smooth projective manifold of dimension $n$. Then any orbifold morphism $f:\bC^n \rightarrow (X/\Delta)$ is degenerate i.e its Jacobian vanishes identically.
\end{theorem}

In \cite{C04} the following orbifold generalization of \cite{KO75} was proved
\begin{theorem}(Campana)
Let $V$ be a connected complex manifold, $Z$ a reduced divisor on $V$ and $U:=V \setminus Z$. Let $\varphi:U \rightarrow X$ be a meromorphic map with $X$ a projective manifold. Let $f:X \rightarrow Y$ be a fibration of general type and assume $\psi:=f\circ \varphi: U \rightarrow Y$ is of maximal rank. Then
\begin{enumerate}
\item $\psi$ extends meromorphically to $V$.
\item for any $m>0$ sufficiently divisible and $s \in H^0(Y,m(K_{Y}+\Delta(f)))$, $\psi^*(s)$ extends to a global holomorphic section of $(\Omega_{V}^p)^{\otimes m}((m-1)Z).$
\end{enumerate}
\end{theorem}

Here we would like to prove the following result which generalizes both \cite{Sa74} and  \cite{C04}
\begin{theorem}\label{KOch}
Let $V$ be a connected complex manifold, $Z$ a reduced divisor on $V$ and $U:=V \setminus Z$. Let $\varphi:U \rightarrow (X/\Delta)$ be an orbifold morphism with $X$ a projective manifold. Let $f:(X/\Delta) \rightarrow Y$ be a fibration of general type and assume $\psi:=f\circ \varphi: U \rightarrow Y$ is of maximal rank. Then
\begin{enumerate}
\item $\psi$ extends meromorphically to $V$.
\item for any $m>0$ sufficiently divisible and $s \in H^0(Y,m(K_{Y}+\Delta(f)))$, $\psi^*(s)$ extends to a global holomorphic section of $(\Omega_{V}^p)^{\otimes m}((m-1)Z).$
\end{enumerate}
\end{theorem}

An immediate corollary which provides new examples of special orbifolds is
\begin{corollary}
Let $(X/\Delta)$ be a smooth orbifold where $X$ is a smooth projective manifold of dimension $n$. Assume there exists a non degenerate orbifold morphism $f:\bC^n \rightarrow (X/\Delta)$. Then $(X/\Delta)$ is special.
\end{corollary}

The proof of the theorem will follow the same lines as \cite{C04} except that we have to take into account the orbifold structure on $X$. We need to recall some definitions from  \cite{C07}.

Let $(X/\Delta)$ be a smooth orbifold. Let $(x_{1},\dots,x_{n})$ be local coordinates such that $\Delta$ has equation $$x_{1}^{(1-\frac{1}{m_{1}})}\dots x_{n}^{(1-\frac{1}{m_{n}})}=0.$$
First, let us recall the definition of the orbifold base of a fibration in this setting.
\begin{definition}
Let $(X/\Delta)$ be an orbifold and $f: (X/\Delta) \rightarrow Y$ a fibration i.e a surjective holomorphic map with connected fibers. Then for every irreducible divisor $D \subset Y$ such that $f^{*}(D)=\sum_{j}m_{j}D_{j}+R$ where $R$ denotes the $f$-exceptional part, we define its multiplicity
$$m(f,\Delta,D)=\inf_{j}\{m_{j}.m_{\Delta}(D_{j}\}.$$
where $m_{\Delta}(D_{j})$ is the multiplicity of $D_{j}$ in $\Delta$.

The orbifold base $(Y/\Delta(f))$ is then defined by
$$\Delta(f)=\sum_{D \subset Y}\left(1-\frac{1}{m(f,\Delta,D)}\right)D.$$
\end{definition}

\begin{definition}
For non negative integers $N,q$, $S^{N,q}(X/\Delta)$ is the locally free subsheaf of $S^{N}\Omega_{X}^q(log\lceil \Delta \rceil)$ generated by the elements $x^{\lceil \frac{k}{m} \rceil}\otimes_{l=1}^{l=N}\frac{dx_{J_{l}}}{x_{J_{l}}}$
where the $J_{l}$ are ordered sets of $q$ elements of $\{1,\dots,n\}$, $x^{\lceil \frac{k}{m} \rceil}=\prod_{j=1}^{n}x_{j}^{\lceil \frac{k_{j}}{m_{j}} \rceil}$ for $k_{j}$ the number of times that $j$ appears in $J_{1},\dots,J_{N}$ and $\frac{dx_{J_{l}}}{x_{J_{l}}}=\wedge_{j \in J_{l}}\frac{dx_{j}}{x_{j}}.$ 
\end{definition}

\begin{definition}
Let $g:(X/\Delta) \rightarrow Y$ be a fibration with $Y$ smooth. $g$ is {\bf neat} realtively to $g':(X'/\Delta') \rightarrow Y'$ if there is a commutative diagramm
$$
\xymatrix{
    (X/\Delta) \ar[r]^w \ar[d]_g& (X'/\Delta') \ar[d]^{g'} \\
    Y \ar[r]_v & Y'
  }
$$
such that
\begin{enumerate}
\item $w$ is an orbifold morphism, $v$ and $w$ are bimeromorphic and $w_*(\Delta)=\Delta'$.
\item Every $g$-exceptional divisor is $w$-exceptional.
\end{enumerate}
\end{definition}

In this context we have (see \cite{C07})
\begin{proposition}\label{neat}
For every fibration  $g':(X'/\Delta') \rightarrow Y'$, there exists $g:(X/\Delta) \rightarrow Y$ neat relatively to $g'$.
\end{proposition}

To deal with the orbifold structure on $X$ we will need the following proposition
\begin{proposition}\label{prop1}
Let $g:(X/\Delta) \rightarrow Y$ and $g':(X'/\Delta') \rightarrow Y'$ be fibrations with $\dim Y =p$ such that $g$ is neat relatively to $g'$.
$$
\xymatrix{
    (X/\Delta) \ar[r]^w \ar[d]_g& (X'/\Delta') \ar[d]^{g'} \\
    Y \ar[r]_v & Y'
  }
$$

Let $m>0$ be an integer such that $m\Delta(g)$ is Cartier. 

Then we have an injection of sheaves $$g^{*}(\mathcal{O}_{Y}(m(K_{Y}+\Delta(g)))) \subset S^{m,p}(X/\Delta).$$
\end{proposition}
\begin{proof}
Let $D \subset X$ be an irreducible divisor not $f$-exceptional and $x_{0} \in D$ a generic point with local coordinates $x=(x_{1},\dots,x_{n})$ such that locally an equation of $D$ is $(x_{1}=0)$. Take local coordinates $(y_{1},\dots,y_{p})$ near $g(x_{0})$ such that  $\Delta(g)$ has equation $$y_{1}^{(1-\frac{1}{m_{1}})}=0$$ and $D$ is mapped locally to $y_{1}=0$. So, $g(x)=(x_{1}^{t_{1}},\dots,x_{p})$ with $t_{1}.m'\geq m_{1}$ for $m'=m_{\Delta}(E).$ 
Let $$\omega=\left(\frac{dy_{1}\wedge\dots\wedge dy_{p}} {y_{1}^{(1-\frac{1}{m_{1}})}}\right)^{\otimes m}$$ be a local generator of $\mathcal{O}_{Y}(m(K_{Y}+\Delta(g)))$. Then 
$$g^{*}(\omega)=x_{1}^{\frac{m t_{1}}{m_{1}}}\dots x_{p}^m\left(\frac{dx_{1}\wedge\dots \wedge dx_{p}} {x_{1}\dots x_{p}}\right)^{\otimes m}$$
up to a non zero constant factor. Since $\frac{m t_{1}}{m_{1}} \geq \lceil \frac{m}{m'} \rceil$, $g^{*}(\omega)$ is a local section of $S^{m,p}(X/\Delta)$. Therefore the injection is true outside $V \cup E(g)$ where $V$ is a subset of $X$ of codimension two or more contained above $\Delta(g)$ and $E(g)$ is the union of $g$-exceptional divisors. But then $w_{*}(g^{*}(\mathcal{O}_{Y}(m(K_{Y}+\Delta(g)))))$ injects into $S^{m,p}(X'/\Delta')$ over $X'$ outside a codimension $2$ or more analytic subset $Z$. Since these sheaves are locally free this injection extends through $Z$ by Hartog extension theorem.
\end{proof}

Now we give the proof of theorem \ref{KOch}
\begin{proof}
We start with some reductions as in \cite{Ko98} and \cite {C04} to which we refer for details. We can reduce to the equidimensional case $\dim U=p$. By localizing we may assume that $V$ is a unit polydisc $D^p$ and $Z$ is a subpolydisc $\{0\}\times D^{p-1}$ so that $U=D^*\times D^{p-1}.$ By proposition \ref{neat} we can assume that $f$ is neat.

Since $K_{Y}+\Delta(g)$ is big, we have
$$m(K_{Y}+\Delta(g))=H+A$$
with $A$ very ample and $H$ effective.
Let $\alpha \in H^0(Y,\mathcal{O}_{Y}(H))$ be a nonzero section and $\delta \in H^0(Y,\mathcal{O}_{Y}(m\Delta))$ a section vanishing exactly on $m\Delta$. Let $\sigma_{0},\dots,\sigma_{N}$ be a basis of $H^0(Y,A)$ and $s_{j}=\alpha\sigma_{j}$ and $t_{j}=\frac{s_{j}}{\delta}$. 
To prove that $\psi$ extends meromorphically to $V$, it suffices to show that the $\psi^*(s_{j})$ extends to meromorphic sections of $mK_{V}$.

Then we define
$$v:=\left(\sum_{j}i^{mp^2}t_{j} \wedge \overline{t_{j}}\right)^\frac{1}{m}$$
which is a meromorphic pseudo-volume form on $Y$.
From proposition \ref{prop1} we see that $f^*(t_{j})$ is a holomorphic section of  $S^{m,p}(X/\Delta)$ and therefore $w=\psi^*(v)$ is a pseudo-volume form on $U$.
In the same way 
$$\tilde{w}=\left(\sum_{j}i^{mp^2}\psi^*(s_{j}) \wedge \overline{\psi^*(s_{j})}\right)^\frac{1}{m}$$
is a pseudo-volume form on $U$.
Then a computation of the Ricci curvature gives
$$Ricci(\tilde{w})=Ricci(w)=\psi^*(-F^*(\Theta))$$
where $F: Y \rightarrow \bP^N$ is the map defined by $A$ and $\Theta$ is the curvature form of the Fubini-Study metric.
Therefore there exists $C>0$ such that $K_{\tilde{w}}\leq -\frac{1}{C}$ and the Schwarz lemma for volume elements imply that 
$$\int_{U} \tilde{w}<\infty.$$
Finally this implies that $\psi^*(s_{j})$ extends to meromorphic sections of $mK_{V}$.

For the second assertion we refer to proposition 8.28 of \cite{C04} 
\end{proof}

\subsection{Orbifold measure hyperbolicity}
\begin{definition}
Let $(X/\Delta)$ be an orbifold with $\Delta=\sum_{i} a_{i}Z_{i}$ and $\Delta_{1}$ the union of all $Z_{i}$ with $a_{i}=1$. 
\begin{enumerate}
\item The orbifold Kobayashi measure $\mu_{(X/\Delta)}$ on $(X/\Delta)$ is the largest measure on $X\setminus \Delta_{1}$ such that
\begin{displaymath}
g^*\mu_{(X/\Delta)} \leqslant \mu_{P}
\end{displaymath}
for every orbifold morphism $g: D^n \rightarrow (X/\Delta)$, where $\mu_{P}$ denotes the Poincar\'e measure on $D^n$ induced by the Poincar\'e volume element $$\kappa_{n}=n!\prod_{j=1}^{n}\frac{4}{(1-|z_{j}|^2)^2}\frac{i}{2}dz_{j}\wedge d\overline{z_{j}}.$$

\item The classical orbifold Kobayashi pseudo-distance $\mu^*_{(X/\Delta)}$ on $(X/\Delta)$ is the largest measure on $X\setminus \Delta_{1}$ such that
\begin{displaymath}
g^*\mu_{(X/\Delta)} \leqslant \mu_{P}
\end{displaymath}
for every classical orbifold morphism $g: D^n \rightarrow (X/\Delta)$.
\end{enumerate}
\end{definition}

\begin{definition}
Let $(X/\Delta)$ be an orbifold. Then $(X/\Delta)$ is said to be (classically) measure hyperbolic if ($\mu^{*}_{(X/\Delta)}(B)>0$) $\mu_{(X/\Delta)}(B)>0$ for every non empty open subset $B \subset  X\setminus \Delta_{1}$.
\end{definition}

\begin{theorem}
Let $(X/\Delta)$ be a projective orbifold of general type. Then $(X/\Delta)$ is measure hyperbolic.
\end{theorem}

\begin{proof}
This is essentially a consequence of \cite{Sa74}. Indeed, there, a volume form $\Psi$ is constructed on $X \setminus supp(\Delta)$ such that $\int_{X \setminus supp(\Delta)} \Psi <\infty$ and for every orbifold morphism $f: D^n \rightarrow (X/\Delta)$, $f^*(\Psi) \leq \kappa_{n}.$ We recall its construction for the convenience of the reader.
Let $A$ be an ample divisor on $X,$ then we can find an effective divisor $D$ in
$$|mK_{(X/\Delta)}-A|$$
for $m$ large enough. Let $\sigma$ be a section vanishing exactly on $D$.
Let $\Delta=\sum_{i} \left(1-\frac{1}{m_{i}}\right)Z_{i}$ with sections $s_{j}$ defining $Z_{j}$.
We choose hermitian metrics $h_{A}$ on $A$, $h_{j}$ on $Z_{j}$. Then we have a volume form on $X \setminus supp(\Delta)$ defined, up to a constant $c$, locally by
$$\Psi=c\frac{||\sigma||_{h_{A}^{-1}}^{\frac{2}{m}}} {\prod_{j=1}^{n} (\ln ||s_{j}||_{h_{j}}^2)^2  |s_{j}|^{2(1-\frac{1}{m_{j}}) }}\prod_{j=1}^{n}\frac{i}{2\pi}dz_{j}\wedge d\overline{z_{j}}.$$
For every orbifold morphism $f: D^n \rightarrow (X/\Delta)$, $f^*(\Psi)$ is a singular volume element such that $(-Ricci(f^*(\Psi)))^n\geq f^*(\Psi)$ in the sense of currents. Therefore, using Ahlfors-Schwarz lemma for volume elements, we obtain that $f^*(\Psi) \leq \kappa_{n}.$ Finally, by definition, $\mu_{(X/\Delta)}$ is greater than the measure induced by $\Psi$ and $(X/\Delta)$ is measure hyperbolic.
\end{proof}

\bigskip
\noindent \texttt{rousseau@math.u-strasbg.fr}

\noindent D\'{e}partement de Math\'{e}matiques,

\noindent IRMA,\newline
Universit\'{e} Louis Pasteur,

\noindent 7, rue Ren\'{e} Descartes,\newline
\noindent 67084 STRASBOURG CEDEX

\noindent FRANCE


\begin{thebibliography}{99999}
\bibitem{Bog77} F.A. Bogomolov, \textit{Families of curves on a surface of general type}, Soviet Math. Dokl. {\bf18} (1977), 1294-1297.

\bibitem{Bog} F.A. Bogomolov, \textit{Holomorphic tensors and vector bundles on projective varieties}, Math. USSR Izvestija {\bf13} (1979), 499-555.

\bibitem{BT} F.A. Bogomolov, Y. Tschinkel, \textit{Special elliptic fibrations}, in Proc. Fano Conf., Torino (2003), ed. A. Conte, arxiv: math.AG/0303044.

\bibitem{Br} M. Brunella, \textit{Courbes enti\`eres et feuilletages holomorphes}, L'Enseignement Math\'ematique {\bf45} (1999), 195-216.

\bibitem{C04} F. Campana, \textit{Orbifolds, special varieties and classification theory}, Ann. Inst. Fourier {\bf 54} (2004), 499-665.

\bibitem{C07} F. Campana, \textit{Orbifoldes sp\'eciales et classification bim\'eromorphe des vari\'et\'es k\"ahl\'eriennes compactes}, arXiv:0705.0737.


\bibitem{CP} F. Campana, M. Paun, \textit{Vari\'et\'es faiblement sp\'eciales \`a courbes enti\`
eres d\'eg\'en\'erees},  Compos. Math.  {\bf143}  (2007),  no. 1, 95--111.

\bibitem{CW} F. Campana, J. Winkelmann, \textit{A Brody theorem for orbifolds}, preprint 2006, arxiv: math/0604571.

\bibitem{ch01} X. Chen, \textit{On Algebraic Hyperbolicity of Log Varieties},
Commun. Contemp. Math.  {\bf 6}  (2004),  no. 4, 513--559. 

\bibitem{De95}  J.-P. Demailly, \textit{Algebraic criteria for Kobayashi
hyperbolic projective varieties and jet differentials}, Proc.
Sympos. Pure Math., vol.62, Amer. Math.Soc., Providence, RI, 1997,
285--360.

\bibitem{DL96}  G. Dethloff, S. Lu, \textit{Logarithmic jet bundles and
applications}, Osaka J. of Math. {\bf38}, 2001, 185-237.

\bibitem{EG} J. El Goul, \textit{Logarithmic jets and hyperbolicity}, Osaka J. Math. {\bf 40} (2003), 469--491.

\bibitem{GG80}  M. Green, P. Griffiths, \textit{Two applications of
algebraic geometry to entire holomorphic mappings}, The Chern Symposium
1979, Proc. Inter. Sympos. Berkeley, CA, 1979, Springer-Verlag, New-York,
1980, 41-74.

\bibitem{Far}  H. M. Farkas, I. Kra, \textit{Riemann Surfaces},
Springer-Verlag, New-York, 1980, second edition.

\bibitem {Jou78} J.-P. Jouanolou, \textit{Hypersurfaces solutions d'une \'equation de Pfaff analytique}, Math. Ann. {\bf 232} (1978), 239-245.

\bibitem{Ko98}  S. Kobayashi, \textit{Hyperbolic complex spaces}, Springer-Verlag, Berlin, 1998.

\bibitem{KO75} S. Kobayashi, T. Ochiai \textit{Meromorphic mappings into compact complex spaces of general type}, Inv. Math., {\bf 31} (1975), 7-16.

\bibitem{McQ0} M. McQuillan, \textit{Diophantine approximations and foliations}, Publ. IHES {\bf 87} (1998), 121-174.

\bibitem{McQ01} M. McQuillan, \textit{Noncommutative Mori theory}, preprint IHES (2000). 

\bibitem{McQ} M. McQuillan, \textit{Bloch hyperbolicity}, preprint IHES (2001).

\bibitem{McQ2} M. McQuillan, \textit{Rational criteria for hyperbolicity}, Book preprint.

\bibitem{Nev} R. Nevanlinna, \textit{Analytic functions}, Berlin-Heidelberg-New York, Springer (1970).

\bibitem{No83} E.I. Nochka, \textit{On the theory of meromorphic functions}, Soviet Math. Dokl. {\bf 27} (2) (1983), 377-381.

\bibitem{PR} G.  Pacienza, E. Rousseau, {\it On the logarithmic Kobayashi conjecture}, 
 J. Reine Angew. Math. {\bf 611}
(2007),  221--235.

\bibitem{Rou03} E. Rousseau \textit{Hyperbolicit\'e du compl\'ementaire d'une courbe : le cas de deux composantes}, CRAS Ser. I 336 (2003), 635-640.

\bibitem{Sa74} F. Sakai, \textit{Degeneracy of holomorphic maps with ramification}, Inv. Math. {\bf26} (1974), 213-229.

\bibitem{Siu} Y.T. Siu, \textit{A proof of the generalized Schwarz lemma using the logarithmic derivative lemma.}  Private communication to J.-P. Demailly, Journal de la SMF (1997).

\bibitem{Vo} P. Vojta, \textit{On the ABC conjecture and diophantine approximation by rational points}, Amer. J. Math.  {\bf122}  (2000),  no. 4, 843--872.

\bibitem{W99} P.-M. Wong, \textit{Nevanlinna theory for holomorphic curves in projective varieties}, preprint (1999).

\end{thebibliography}
\end{document}